\documentclass[oneside,a4paper,12pt,reqno]{amsart}
\newcommand{\tor}{\operatorname{Tor}}
\newcommand{\asc}{\operatorname{Asc}}

\newtheorem{theorem}{Theorem}[section]
\newtheorem{proposition}[theorem]{Proposition}
\newtheorem{lemma}[theorem]{Lemma}

\theoremstyle{definition}
\newtheorem{example}[theorem]{Example}
\def\eqref#1{(\ref{#1})}

\begin{document}

\bibliographystyle{amsplain}

\title[Topological rigidity on connected groups]{Finite
entropy characterizes
topological rigidity on connected groups}

\author{Siddhartha Bhattacharya}
\address{SB: School of Mathematics,
Tata Institute of Fundamental Research, Bombay 400005, India}
\email{siddhart@math.tifr.res.in}
\author{Thomas Ward}

\address{TW: School of Mathematics, University of
East Anglia, Norwich NR4 7TJ, England} \email{t.ward@uea.ac.uk}
\thanks{The first author gratefully acknowledges the
support of London Mathematical Society grant~4714 and the
hospitality of the University of East Anglia.}
\begin{abstract}
Let $\mathsf X_1$, $\mathsf X_2$ be mixing connected
algebraic dynamical systems with the
Descending Chain Condition. We show that every
equivariant continuous map $\mathsf X_1\to\mathsf X_2$
is affine (that is, $\mathsf X_2$ is {\sl topologically
rigid}) if and only if the system $\mathsf X_2$ has
finite topological entropy.
\end{abstract}
\maketitle
%
%==========================
%  Section 1
%==========================
%
\section{Introduction}
An {\it algebraic ${\mathbb{Z}}^{d}$-action} $\alpha $
on a compact abelian group $X$ is a
homomorphism $\alpha:{\mathbf n}\mapsto\alpha({\mathbf n})$
from ${\mathbb{Z}}^{d}$ to the group ${\rm Aut}(X)$ of
continuous automorphisms
of $X$. Compact groups are
assumed to be metrizable throughout and
are written multiplicatively; $e$ is used to
denote the identity element of any group.
Write $\mathsf X=(X,\alpha)$ for such an
algebraic dynamical system, and call the
system $\mathsf X$ connected, mixing and so on if
$X$ is connected, $\alpha$ is mixing, and so on.

Any algebraic system $\mathsf X$ preserves
$\lambda_{X}$, the Haar measure on $X$.
The system $\mathsf X$ is {\sl mixing} if
$$
\lim\limits_{\mathbf n\to\infty}
\lambda_{X}(A_{1}\cap\alpha({\mathbf n})(A_{2}))
=\lambda_{X}(A_{1})\cdot\lambda_{X}(A_{2})
$$
for all measurable sets $A_1,A_2\subset X.$

A map $\phi:\mathsf X_1\to\mathsf X_2$ between
algebraic dynamical systems is {\sl equivariant}
if $\phi\circ\alpha_1(\mathbf n)=\alpha_2(\mathbf n)\circ\phi$
for all $\mathbf n\in\mathbb Z^d$, and is {\sl affine}
if there is a continuous group homomorphism
$\psi:X_1\to X_2$ and an element $y\in X_2$ with
$\phi(x)=\psi(x)\cdot y$.

Topological (respectively, measurable) rigidity is a property of the
target system $\mathsf X_2$ that forces an equivariant
continuous (resp. measurable) map to coincide everywhere
(resp. almost everywhere) with an affine map.

For $d\ge 1$, denote by
$R_d={{\mathbb{Z}}}[u_1^{\pm 1},\dots,u_d^{\pm 1}]$
the ring of Laurent polynomials with integral coefficients in
$d$ commuting variables $u_1,\ldots,u_d$.
An element $f$ of $R_d$ is written
        \begin{equation*}
f(\mathbf u)=\sum_{{\mathbf n}\in{{\mathbb{Z}}}^d}
f_{\mathbf n}{\mathbf u}^{\mathbf n}
        \end{equation*}
with ${\mathbf u}^{\bf n}=u_1^{n_1}\cdots u_d^{n_d}$,
$f_{\mathbf n}\in\mathbb{Z}$ for
all ${\mathbf n}=(n_1,\dots ,n_d)\in\mathbb{Z}^d$, and $f_{\mathbf n}=0$
for all but finitely many ${\mathbf n}\in\mathbb{Z}^d$.

If $\mathsf X=(X,\alpha)$ is an algebraic $\mathbb Z^d$-action
on a compact abelian group $X$, then the countable
dual group $M=\widehat{X}$ is a
module over the ring $R_d$ under the operation
$$
f\cdot a=\sum_{{\mathbf n}\in{{\mathbb{Z}}}^d}
f_{\mathbf n}{\widehat{\alpha}}
({\mathbf n})(a)
$$
for $f\in R_d$ and $a\in M$.
The module $M$ is called the \emph{dual module} of $\mathsf X$.
Conversely, a countable module $M$ over $R_d$ determines an
algebraic $\mathbb Z^d$-action $\mathsf X_M=(X_M,\alpha_M)$
by setting
$${\widehat{\alpha _M}}({\mathbf n})(a)={\mathbf u}^{\mathbf n}\cdot a $$
for every ${\mathbf n}\in {{\mathbb{Z}}}^{d}$ and $a\in M$.

An algebraic $\mathbb Z^d$-action $\mathsf X$ is {\sl Noetherian}
if the dual module is Noetherian.
The following properties are equivalent.
\begin{itemize}
\item $\mathsf X_M$ is Noetherian.
\item $M$ is finitely-generated over $R_d$ (this is equivalent
to $M$ being Noetherian since $R_d$ is itself Noetherian).
\item Any descending chain of
closed $\alpha_M$-invariant subgroups of $X_M$
stabilizes (the Descending Chain Condition;
see~\cite{MR91g:22008}).
\end{itemize}

The topological entropy of the system $\mathsf X_M$ is
defined and computed in terms of the module $M$ in~\cite{MR92j:22013}.

Rigidity properties of algebraic ${\mathbb{Z}}^{d}$-actions
have been
studied by several authors.
Measurable equivariant
maps between mixing zero-entropy
algebraic ${\mathbb{Z}}^{d}$-actions exhibit strong regularity
properties (see~\cite{B2}, \cite{BS}, \cite{KKS} and~\cite{MR2001j:37004}).
For a $\mathbb Z$-action generated by an
automorphism $\theta$
on a connected finite-dimensional compact abelian group, it
is known that the topological centralizer of the action admits
non-affine maps if and only if $\theta$ is not
ergodic (cf.~\cite{MR33:1402}, \cite{MR27:5858}
and~\cite{MR39:1631}). Ergodic automorphisms of
infinite-dimensional groups may have non-affine
maps in their centralizers (see Example~\ref{infinitedim}).
In~\cite{MR2000m:37017}
it is shown that for any expansive connected algebraic
${\mathbb{Z}}^{d}$-action $\mathsf X$, the
topological centralizer of $\alpha$ consists of affine maps
(expansiveness is a condition that implies the
Descending Chain Condition; for $d=1$ it
forces the compact group $X$ to be finite-dimensional).

In this note we prove the following result, which
characterizes a form of topological rigidity in terms of
topological entropy.
\begin{theorem}\label{rigidity}
Let $\mathsf X_{1}$, $\mathsf X_{2}$ be connected
mixing Noetherian algebraic ${\mathbb{Z}}^{d}$-actions.
Then the following properties are equivalent.
\begin{enumerate}
\item Every equivariant continuous map
$\mathsf X_{1}\to\mathsf X_{2}$ is an affine map.
\item The system $\mathsf X_{2}$ has finite topological entropy.
\end{enumerate}
\end{theorem}
For $d\ge2$, this result applies to situations where
the underlying group $X_2$ is infinite-dimensional,
so the lifting techniques of~\cite{MR33:1402} and~\cite{MR39:1631} cannot
be applied directly.
Write $\mathbb T\subset\mathbb C$ for
the multiplicative unit circle.

\begin{example}\label{ledrappier}
To illustrate Theorem~\ref{rigidity},
consider the $\mathbb Z^2$-action $\mathsf X$ where
$X\subset\mathbb T^{\mathbb Z^2}$ is the closed subgroup
consisting of all $x\in {{\mathbb{T}}}^{{\mathbb{Z}}^{2}}$
with
$$
x(m+1,n)\cdot x(m,n)\cdot x(m,n+1) = 1\mbox{ for all }
m,n\in\mathbb Z,$$
and $\alpha$ is the shift action of ${\mathbb{Z}}^{2}$ on $X$.
The system $\mathsf X$ is mixing and has finite
entropy. It follows that every
continuous equivariant map
from $\mathsf X$ to itself
is an affine map.
In contrast, the {\sl measurable} centraliser of $\mathsf X$
contains many non-affine maps, since $\mathsf X$ is
measurably isomorphic to a
$\mathbb Z^2$ Bernoulli shift (see~\cite{MR96d:22004}, \cite{MR93h:28032}).
\end{example}

\begin{example}\label{infinitedim}
For the case of a single automorphism,
the compact group being finite-dimensional
forces the entropy to be finite.
Ergodic automorphisms of infinite-dimensional
groups are not topologically rigid in general.

For example, the shift
automorphism of $X=\mathbb T^{\mathbb Z}$
defines
an ergodic $\mathbb Z$-action of infinite
entropy that is not
topologically rigid: if $f:\mathbb T\to\mathbb T$
is any map, then the shift map
commutes with the map $\phi:X\to X$ defined
by $(\phi(x))_k=f(x_k)$. The module corresponding
to this action is a Noetherian $R_1$-module.

On the other hand, an ergodic automorphism of
$\mathbb T^{\mathbb Z}$ that splits into
a direct product of automorphisms of
finite-dimensional tori {\sl is} topologically
rigid. The module corresponding to this action
is not Noetherian.
It is not know whether such an action
can have finite topological entropy (see~\cite{MR49:10856}).
\end{example}

The next example again shows that Theorem~\ref{rigidity}
does not hold for non-Noetherian actions.

\begin{example}\label{last}
Let $F_{d}$ denote the field
of fractions of $R_{d}$, considered as a $R_d$-module.
Let $\mathsf X_1$ denote the algebraic
$\mathbb Z^d$-action
corresponding to $F_d$.
Notice that $F_d$ is torsion-free as a
$R_d$-module, and $\mathsf X_1$
has infinite entropy. For any
$\mathbf n\in\mathbb Z^d$, multiplication by
${\mathbf u}^{\mathbf n}-1 $ is an automorphism of $F_d$.
By duality,
the map $x\mapsto\alpha_1(\mathbf n)(x)-x$
is a continuous automorphism
of $X_1$ for any $\mathbf n\in\mathbb Z^d$.
In particular, $\mathsf X_1$ does not have any
non-trivial periodic orbits. Now let
$\mathsf X_2$ be any mixing connected
algebraic ${\mathbb{Z}}^{d}$-action with
a dense set of periodic orbits (any Noetherian
system has this property). Since continuous
equivariant maps take periodic
orbits to periodic orbits, it follows that any continuous equivariant
map from $\mathsf X_2$ to $\mathsf X_1$ is trivial.
\end{example}

Example~\ref{last} is similar in spirit to
a remark of Comfort (see~\cite{MR36:2734}): there are
no non-trivial homomorphisms $\mathbb T\to\widehat{\mathbb Q}$
since torsion elements are dense in $\mathbb T$ but
absent in $\widehat{\mathbb Q}$.

%
%==========================
%  Section 2
%==========================
%
\section{Algebraic $\mathbb Z^d$-actions}

In this section basic results and terminology on
algebraic $\mathbb Z^d$-actions is collected.
A prime ideal
$\mathfrak p\subset R_d$ is {\sl associated with} the
$R_d$-module $M$
if there exists $m\in M$ with
${\mathfrak{p}}=\{ f\in R_{d}\ |\ f\cdot m = 0\}$.
The set of prime ideals associated with $M$
is denoted $\asc(M)$.
If $M$ is Noetherian, then $\asc(M)$ is finite.
The torsion
submodule of $M$ is defined by
$$\tor(M)=\{m\in M\mid r\cdot m=0\mbox{ for some non-zero } r\in R_d\}.$$
A module $M$ is said to be a {\sl torsion module}
if $\tor(M)=M$.

The following result taken from~\cite[Theorem~6.5]{MR97c:28041}
characterizes mixing in algebraic terms.

\begin{lemma}
The algebraic $\mathbb Z^d$-action $\mathsf X_M$ is mixing if and
only if for every $\mathfrak{p}\in\asc(M)$ and for
every non-zero $\mathbf n\in\mathbb Z^d$,
the polynomial ${\mathbf u}^{\mathbf n}-1$ does not lie in
$\mathfrak p$.
\end{lemma}

An algebraic $\mathbb Z^d$-action $\mathsf X_2$ is
an {\it algebraic factor} of
$\mathsf X_1$ if there is a surjective
continuous equivariant homomorphism $\phi:X_1\to X_2$.

The next lemma shows that if the module corresponding to an
algebraic $\mathbb Z^d$-action is Noetherian,
then infinite topological entropy can only
be created by the presence of (a factor of)
a full shift with infinite alphabet.

\begin{lemma}\label{ent}
For a Noetherian system $\mathsf X_M$
the following conditions are
equivalent.
\begin{enumerate}
\item $X_M$ does not admit a non-trivial closed $\alpha_M$-invariant subgroup
$H$ with the property that
the restriction of $\alpha_M$ to $H$ is an
algebraic factor of the shift action of $\mathbb Z^d$ on
$(\mathbb T^n)^{\mathbb Z^{d}}$ for some $n>0$.
\item $M$ is a torsion module.
\item $\mathsf X_M$ has finite topological entropy.
\end{enumerate}
\end{lemma}

\begin{proof}
$(1)\implies(2).$ Suppose that $M$ is not a torsion module
and $N=M/\tor(M)$.
Then $N$ is a non-zero torsion-free $R_d$-module.

We claim that $N$
is isomorphic to a submodule of the free module $R_d^n$
of rank $n$ for some $n\ge1$.
Let $N_{0}$ denote the localisation of $N$ at the prime ideal $\{0\}$.
Since $N$ is Noetherian, $N_{0}$ is a finite-dimensional
vector space over $F=\mathbb Z(u_1^{\pm1},\dots,u_d^{\pm1})$,
the quotient field of $R_d$. Let
$B=\{b_1,\dots,b_n\}$ be any $F$-basis of $N_{0}$.
The map $m\mapsto m/1$ embeds $N$ as a submodule of
$N_0$.
Choose a finite $R_d$-generating set $A$ of $N$, and
an element $p\in R_d$ with the property that
$p\cdot a$ lies in the $R_d$-submodule generated
by $B$ for all $a\in A$.
The submodule generated by
$\{b_{1}/p,\dots,b_{n}/p\}$ contains $N$, and is a free
$R_{d}$-module of rank $n$.
This proves the claim.

The system $\mathsf X_N$ is therefore an
algebraic factor of the shift action on
$(\mathbb T^n)^{\mathbb Z^d}$.
Since $N$ is a quotient of $M$, by duality there exists
a closed $\alpha$-invariant subgroup
$H\subset X$ such that the restriction of $\alpha$ to $H$ is
conjugate to $\mathsf X_N$.
\smallskip

$(2)\implies(1).$
Let $H\subset X$ be a closed $\alpha$-invariant subgroup such that
the restriction of $\alpha$ to $H$ is an
algebraic factor of the shift action on
$(\mathbb T^n)^{\mathbb Z^d}$ for some $n>0$.
The dual module of the shift action
on $({\mathbb{T}}^{n})^{{\mathbb{Z}}^{d}}$
is isomorphic to the direct sum of $n$ copies of $R_d$,
which implies that $\widehat{H}$ is a torsion-free $R_d$-module.
On the other hand,
$\widehat{H}$ is a quotient of the $R_d$-module $M$,
which is a torsion module by assumption (2).
Hence $\widehat{H}=\{0\}$ so $H$ is trivial.

$(2)\implies(3).$
Let $\{m_{1},\dots,m_{k}\}$ generate
$M$ as an $R_d$-module. For $j=1,\dots,k$ let
$I_j\subset R_d$ denote the ideal defined by
$$
I_j=\{p\in R_d\mid p\cdot m_j=0\}.
$$
Since $M$ is a torsion module, each $I_j$ is non-zero.
For $j=1,\dots,k$ let $M_j$ denote the
$R_d$-module $R_d/I_j$.
Since each $I_j$ is non-zero,
$\mathsf X_{M_j}$ has finite entropy by~\cite[Theorem~3.1]{MR92j:22013}.
Let
$$
M^{\prime}=M_1\oplus\dots\oplus M_k.
$$
Since each $\mathsf X_{M_j}$ has finite entropy,
$\mathsf X_{M^{\prime}}$ also has
finite entropy.
The map $(r_1,\dots,r_k)\mapsto
r_1m_1+\dots+r_km_k$ expresses $M$ as a quotient of
$M^{\prime}$. The dual of
this map embeds $\mathsf X_M$ as a
sub-action of $\mathsf X_{M^{\prime}}$, so in particular
$\mathsf X_M$ has finite entropy.

$(3)\implies(2).$
If $M$ is not a torsion module, then it contains
$R_d$ as a submodule. By duality,
the shift action of ${\mathbb{Z}}^{d}$ on
${\mathbb{T}}^{{\mathbb{Z}}^{d}}$ is therefore an algebraic
factor of $\mathsf X_M$. Since the former action
has infinite entropy, $\mathsf X_M$ has infinite
topological entropy.
\end{proof}

\section{van Kampen's theorem}

In the proof of Theorem~\ref{rigidity} the following
structure theorem of van Kampen~\cite{VK} will be used
in place of the lifting of toral maps. This
result splits continuous maps into a `linear'
part (a character) and a `non-linear' part in
a unique way. It is also used in this connection by
Walters~\cite{MR39:1631}.
We include a short proof
for the convenience of the reader.
\begin{theorem}\label{VK}
Let $X$ be a compact connected
abelian group and let $f:X\rightarrow\mathbb{T}$ be a continuous map
with $f(e)=1$.
Then there exist a character $\phi\in{\widehat X}$
and a continuous map
$S(f):X\rightarrow\mathbb{R}$ such that
\begin{equation}\label{vankampen}
S(f)(e)=0,\quad f(x)=\phi(x)\cdot e^{2\pi iS(f)(x)}\mbox{ for all }
x\in X.
\end{equation}
Moreover, $\phi$ and $S(f)$ are uniquely defined by~{\rm\eqref{vankampen}}.
\end{theorem}

\begin{proof}
Write $\exp(t)=e^{2\pi it}$.
We first prove that~\eqref{vankampen} determines
$\phi$ and $S(f)$ uniquely.
Suppose that
$$f=\phi_{1}\cdot(\exp\circ h_{1})=\phi_{2}\cdot(\exp\circ h_{2}).$$
Then $\phi\cdot(\exp\circ h)=1$, where
$\phi=\phi_{1}-\phi_{2}$ and $h=h_{1}-h_{2}$.
Define $p:X\times X\rightarrow\mathbb R$ by
$$p(x,y)=h(x+y)-h(x)-h(y).$$
Since $\phi$ and $\exp$ are homomorphisms
and $\exp\circ h={\overline{\phi}}$,
it follows that $\exp\circ p$ is trivial, so the image of $p$
is contained in $\mathbb Z$. Since $X$ is connected, this implies that
$p$ is identically zero, so $h$ is a continuous homomorphism from
$X$ to $\mathbb R$. This forces $h$ to be identically zero
(since $\mathbb R$ has no non-trivial compact subgroups), so
$h_{1}=h_{2}$ and $\phi_{1}=\phi_{2}$.

Turning to the existence of $\phi$ and $S(f)$,
we first consider two special cases.
\begin{enumerate}
\item If $X=\mathbb T^n$ for some $n$, then
let $\pi_{1}(f):\mathbb Z^n\rightarrow\mathbb Z$
be the induced map between fundamental groups.
Choose $\phi\in\widehat X$ with the property
that $\pi_{1}(f)=\pi_{1}(\phi)$.
Since the map $f\cdot\overline{\phi}$ is null-homotopic
and $\exp:\mathbb R\rightarrow\mathbb T$
is a covering map,
there exists a unique continuous map $h:X\rightarrow\mathbb R$
such that $h(e)=0$ and $f\cdot\overline\phi=\exp\circ h$.
\item If $f(X)\subset V=\exp(-1/4,1/4)$ then, since
the restriction of $\exp$
to $(-1/4,1/4)$ is one-to-one, there exists a unique map
$$p:V\rightarrow(-1/4,1/4)$$
such that $\exp\circ p$ is the identity on $V$.
Thus $f=\exp\circ p\circ f$.
\end{enumerate}
To prove the general case, choose a translation-invariant metric $\rho$
on $X$. Since
$f$ is uniformly continuous, there exists $\delta>0$ such that
$$
\rho(x,y)<\delta\implies f(x){\overline{f(y)}}\in V.
$$
For any finite set $A\subset\widehat{X}$,
let
$$
X(A)=\{x\in X\mid\phi(x)=1\mbox{ for all }\phi\in A\}.
$$
Then $X(A)$ is a closed subgroup of $X$ since it is
a finite intersection of kernels of characters.
Moreover, the family of sets
$$
\{X(A)\mid A\subset\widehat{X}\mbox{ is finite}\}$$
has the finite intersection property,
and their intersection is $\{e\}$ since the
character group separates points.

Let $Y\subset X$ denote the open subset defined by
$$Y=\{x\in X\mid\rho(e,x)<\delta/2\}.$$
Since $Y^{c}=X\backslash Y$ is compact,
there exists a finite set $B\subset\widehat{X}$
such that $X(B) \subset Y$.
For any $x\in X$, let
$$A_x=\{f(xz)\mid z\in X(B)\}
\subset\mathbb{T}.
$$
Since the diameter of $A_x$ is strictly less than $1/2$,
there exists a unique closed ball of minimal radius,
$B_x\subset\mathbb T$, that
contains $A_x$.
Let $r(x)$ be the center of $B_{x}$.
The map $x\mapsto r(x)$ is well-defined, continuous,
and invariant under translations by elements of $X(B)$.
Since $X$ is compact and connected, $\widehat{X}$ is
torsion-free. On the other hand,
the annihilator $X(B)^{\perp}$ is
finitely-generated, so
$X/X(B)$ is isomorphic to a torus. It follows that
$r=\phi\cdot(\exp\circ p)$ for some $\phi\in\widehat{X}$
and some $p:X\rightarrow\mathbb R$
with $p(e)=0$.
Since $f\cdot\overline{r}(X)\subset V$,
$f\cdot\overline{r}=\exp\circ q$
for some $q:X\rightarrow\mathbb R$ with $q(e)=0$.
Now $f=r\cdot(f\cdot\overline{r})=\phi\cdot(\exp\circ(p+q))$.
\end{proof}
%
%==============================
%   SECTION 3
%==============================
%
\section{Rigidity of equivariant maps}

For any algebraic $\mathbb Z^d$-action
$\mathsf X=(X,\alpha)$ and for any locally compact
abelian group $A$, denote by $A^{X}$ the
group of all continuous maps
$$
h:X\rightarrow A,\quad h(e)=e,
$$
equipped with point-wise multiplication.
The action $\alpha$ induces the structure of an $R_d$-module
on $A^{X}$ by defining
$$
p\cdot h(x)=\sum_{\mathbf n\in\mathbb Z^d}
p(\mathbf n)\cdot h\circ\alpha(\mathbf n)(x).
$$
A key observation is that
$\widehat{X}$ can be regarded as a submodule
of $\mathbb T^X$ with this structure.

\begin{proposition}\label{ES}
Let $\mathsf X=(X,\alpha)$ be a connected $\mathbb Z^d$-action.
Then $\mathbb R^X$ and $\mathbb T^X/\widehat{X}$
are isomorphic as $R_d$-modules.
\end{proposition}

\begin{proof} The correspondence $f\mapsto S(f)$ from
Theorem~\ref{VK} induces a map $S$ from
$\mathbb T^X$ to $\mathbb R^X$.
If $f_{1}, f_{2}$ are elements of $\mathbb T^X$ then,
by the uniqueness part of Theorem~\ref{VK},
$S(f_1\overline{f_2})=S(f_1)-S(f_2)$
so $S$ is a group homomorphism.
Similarly, if $\theta$ is a continuous endomorphism of $G$,
then $S(f\circ\theta)=S(f)\circ\theta.$
Hence
$S :\mathbb T^X\rightarrow\mathbb R^X$
is an $R_{d}$-module homomorphism.
For any $f$ in $\mathbb R^X$,
$S(e^{2\pi i f})=f$, so the map $S$ is surjective.
Since $\ker(S)=\widehat{X}$, the statement follows.
\end{proof}

If $f$ and $g$ are functions from $\mathbb Z^d$ to
$\mathbb{C}$ and $g$ has finite support,
the convolution $f\ast g:\mathbb Z^d\rightarrow\mathbb C$
is given by
$$
f\ast g(\mathbf i)=\sum_{\mathbf j\in\mathbb Z^d}
f(\mathbf i-\mathbf j)\cdot g(\mathbf j).$$
Write $L^2(\mathbb Z^d)$ for the set of all square-integrable functions
$\mathbb Z^d\to\mathbb C$ (with respect to the
counting measure on $\mathbb Z^d$).

In addition to van Kampen's theorem, a simple
version of the $L^2$ zero-divisor problem is needed
(see~\cite{MR90k:16011} for an overview).

\begin{proposition}\label{zdc}
If $f\in L^2(\mathbb Z^d)$ has
$f\ast g=0$ for some non-zero function
$g :\mathbb Z^d\rightarrow\mathbb C$ with finite support,
then $f$ is identically zero.
\end{proposition}

\begin{proof}
Since the support of $g$ is finite, there exists
$\mathbf n\in\mathbb Z^d$ such that the support of
$g\ast\delta_{\mathbf n}$ is contained in $\mathbb N^d$.
Replacing $g$ by $g\ast\delta_{\mathbf n}$ if necessary,
we may assume that the support of $g$ is contained in
$\mathbb N^d$.
Let $\widehat{f},\widehat{g}\in L^2(\mathbb T^d)$ denote
the Fourier transforms of $f$ and $g$ respectively.
By the choice of $g$, $\widehat{g}=p|_{\mathbb T^d}$ for
some non-zero polynomial $p(\mathbf z)\in\mathbb C[z_1,\dots,z_d]$.
Define
$V(p)\subset\mathbb T^d$ by
$$
V(p)=\{x\in\mathbb T^d\mid p(x)=0\}.
$$
We claim that $\lambda_{d}(V(p))=0$
for any non-zero $p$,
where $\lambda_d$ is Haar measure
on $\mathbb T^d$. This may be proved by induction on $d$.
If $d=1$, then $V(p)\subset\mathbb T$ is finite since
every non-zero polynomial has only finitely many roots.
If $d>1$, choose polynomials
$p_0,\dots,p_k$ in $\mathbb C[z_1,\dots,z_{d-1}]$
such that
$$
p(z_1,\dots,z_d)=\sum_{i=0}^{k}p_i
(z_1,\dots,z_{d-1})z_d^i.
$$
Since $p$ is non-zero, $p_{i}$ is non-zero for some $i$.
By the inductive hypothesis, $\lambda_{d-1}(V(p_i))=0$.
If $(z_1,\dots,z_{d-1})$ lies in $\mathbb T^{d-1}\backslash V(p_i)$,
then the map $z\mapsto p(z_1,\dots,z_{d-1},z)$ is a
non-zero polynomial in $\mathbb C[z]$. This
implies that
for any $(z_1,\dots,z_{d-1})\in\mathbb T^{d-1}\backslash V(p_{i})$,
the set
$$\{z\in\mathbb T\mid(z_1,\dots,z_{d-1},z)\in V(p)\}$$
is finite.
By Fubini's theorem,
$$
\lambda_d(V(p))=\int_{\mathbb T^{d-1}}\int_{\mathbb T^{\vphantom{d-1}}}
\mathbb{I}_{V(p)}\mbox{d}\lambda_{1} d\lambda_{d-1}=0,$$
where $\mathbb{I}$ is the indicator function,
which proves the claim.
Since $\widehat{g}=p|_{\mathbb T^d}$ and
$\widehat{f}\cdot\widehat{g}=\widehat{f\ast g}=0$,
it follows that $\widehat{f}=0$ almost everywhere, so $f=0$.
\end{proof}

\begin{lemma}\label{Mix}
If $\mathsf X=(X,\alpha)$
is a mixing connected algebraic ${\mathbb{Z}}^{d}$-action,
then
$\tor(\mathbb T^X)\subset\widehat{X}$.
\end{lemma}

\begin{proof}
Let $M$ denote the set of all square-integrable
functions
$$
h:\widehat{X}\to\mathbb{C}, \quad h(e)=0.
$$
Defining
$$
p\cdot h(\phi)=\sum_{\mathbf i\in\mathbb Z^d}
p(\mathbf i)h(\phi\circ\alpha(\mathbf i))
$$
for $p\in R_d$ gives $M$ the structure of an $R_d$-module.

We claim first that $M$ is torsion-free.
Let $h$ be an element of $\tor(M)$; for any non-trivial
$\chi\in\widehat{X}$ define a function
$h_{\chi}:\mathbb Z^d\rightarrow\mathbb C$ by
$$
h_{\chi}(\mathbf i)=h(\chi\circ\alpha(\mathbf i)).
$$
Since $\alpha$ is mixing, the map
$\mathbf i\mapsto\chi\circ\alpha(\mathbf i)$ is one-to-one.
Hence
$$
\sum_{\mathbf i\in\mathbb Z^d}|h_{\chi}(\mathbf i)|^2
\ \le \ \sum_{\chi\in\widehat{X}}|h(\chi)|^2<\infty.
$$
This shows that $h_{\chi}\in L^2(\mathbb Z^d)$ for
all $\chi\in\widehat{X}$.
Note that $L^2(\mathbb Z^d)$ itself is
an $R_{d}$-module with respect to the
multiplication $p\cdot h=p\ast h$.
Furthermore, the map
$h \mapsto h_{\chi}$ is an $R_{d}$-module homomorphism
from $M$ to $L^2(\mathbb Z^d)$.
Since $\tor(L^2(\mathbb Z^d))=\{0\}$ by
Proposition~\ref{zdc}, we conclude that $h_{\chi}=0$ for
all $\chi$, so $h=0$. This proves that $M$ is torsion-free.

For $f$ in $\mathbb R^X$, let
$\widehat{f}\in M$ denote the Fourier transform of $f$. Since
$\widehat{f\circ \theta}(\phi)=\widehat{f}(\phi\circ\theta)$
for any continuous endomorphism $\theta$ of $X$,
the map $f\mapsto\widehat{f}$ is an
$R_{d}$-module homomorphism
from $\mathbb R^X$ to $M$.
By the Fourier inversion
theorem this map is injective. Since $\tor(M)=\{0\}$,
this implies
that $\tor(\mathbb R^X)=\{0\}$. Proposition~\ref{ES}
then shows that
$\tor(\mathbb T^X)\subset\widehat{X}$.
\end{proof}

We are now ready to prove Theorem~\ref{rigidity}.

\begin{proof}
Suppose that $\mathsf X_2$ has finite entropy, and let
$f$ be an equivariant continuous
map $\mathsf X_1\to\mathsf X_2$. Define
$f_0:X_1\rightarrow X_2$ by
\begin{equation}\label{fzero}
f_0(x)=f(x)-f(e).
\end{equation}
Since $f$ is equivariant, so is $f_0$.

Fix an arbitrary character $\phi\in\widehat{X_2}$.
By Lemma~\ref{ent}, $\widehat{X_2}$ is a torsion
module, so $\phi$ lies in  the torsion submodule of
$\mathbb T^{X_{2}}$.
Since $f_0$ is equivariant and $f_0(e)=1$,
the map $h\mapsto h\circ f_{0}$ is an $R_{d}$-module
homomorphism $\mathbb T^{X_{2}}\to\mathbb T^{X_{1}}$.
Hence $\phi\circ f_0$ is an element
of the torsion submodule of
$\mathbb T^{X_{1}}$. By Lemma~\ref{Mix},
$\phi\circ f_{0}$ lies in $\widehat{X_1}$. Since the
initial choice of $\phi$ was arbitrary, this shows that
$\phi\mapsto\phi\circ f_{0}$ is a group homomorphism
from $\widehat{X_2}$ to $\widehat{X_1}$. By duality,
there exists a continuous homomorphism $\theta:X_1\to X_2$
such that
$\phi\circ f_0=\phi\circ\theta$
for all $\phi\in\widehat{X_2}.$
Since characters separate points, this implies
that $f_0=\theta$.
Hence $f=f(e)+f_0$ is an affine map.

If $\mathsf X_2$ has infinite entropy,
then by Lemma~\ref{ent} there exists a non-trivial
closed $\alpha_2$-invariant subgroup
$H\subset X_{2}$ with the property that
the restriction of $\alpha_2$ to $H$ is an
algebraic factor of the shift action on
$(\mathbb T^n)^{\mathbb Z^d}$ for some $n>0$.
Let $K\subset(\mathbb T^n)^{\mathbb Z^d}$ be a
proper, closed shift-invariant subgroup such that
the restriction of $\alpha_{2}$ to $H$ is algebraically
conjugate to the shift action of ${\mathbb{Z}}^{d}$ on
$({\mathbb{T}}^{n})^{{\mathbb{Z}}^{d}}/K$. Since any non-affine equivariant
map from $X_{1}$ to $H$ gives rise to a
non-affine equivariant
map $\mathsf X_1\to\mathsf X_2$,
without loss of generality we may assume that
$X_2=(\mathbb T^n)^{\mathbb Z^d}/K$, and
$\alpha_2$ is the shift action.

For any continuous map $q:X_1\rightarrow\mathbb R^n$
define a map
$$
\sigma(q):X_1\rightarrow(\mathbb T^n)^{\mathbb Z^d}
$$
by
$$
\sigma(q)(x)(\mathbf n)=\exp\circ q\circ \alpha_1(\mathbf n)(x).
$$
Let $\pi:(\mathbb T^n)^{\mathbb Z^d}\to
(\mathbb T^n)^{\mathbb Z^d}/K$ denote the
projection map. For any
$$
q:X_1\rightarrow\mathbb R^n,
$$
$\pi\circ\sigma(q)$ is a continuous equivariant map
from $\mathsf X_1$ to $\mathsf X_2$.
We claim that there exists a non-zero
continuous equivariant map
from $\mathsf X_1$ to $\mathsf X_2$ of the form
$\pi\circ\sigma(q)$
for some continuous map $q:X_1\rightarrow\mathbb R^n$
with $q(e)=0$.

For any finite set $F\subset\mathbb Z^d$, let
$\Pi_F$ denote the projection map
$$\Pi_F:(\mathbb T^n)^{\mathbb Z^d}\to
(\mathbb T^n)^F.$$
Since $K$ is a proper closed subgroup of
$(\mathbb T^n)^{\mathbb Z^d}$, there exists
a finite set $F\subset\mathbb Z^d$, and a point
$\mathbf x\in(\mathbb T^n)^F$, such that $\mathbf x$ does not
lie in the image of $\Pi_F$.
Since $\mathsf X_1$ is mixing,
for any $\mathbf i\neq\mathbf j\in\mathbb Z^d$,
the kernel of $\alpha_1(\mathbf i)-\alpha_1(\mathbf j)$
is a proper closed subgroup of $X_1$.
In particular, there exists $y\in X_1$ such that $y\ne e$,
and
$$
\alpha_1(\mathbf i)(y)\ne\alpha_1(\mathbf j)(y)
\mbox{ for any }\mathbf i,\mathbf j\in F.
$$
Choose $\mathbf z\in(\mathbb R^n)^F$ such
that $\exp(\mathbf z(\mathbf i))=\mathbf x(\mathbf i)$ for
all $\mathbf i\in F$. Let
$$q:X_1\to\mathbb R^n$$
be any continuous map with
$q(e)=0$ and $q\circ\alpha_1(\mathbf i)(y)=\mathbf z(\mathbf i)$
for all $\mathbf i\in F$.
Since
$\pi\circ\sigma(q)(y)$ does not lie in $K$, this proves the claim.

Now let $q:X_1\rightarrow\mathbb R^n$
be any continuous map such that $q(e)=0$,
and $\pi\circ\sigma(q):X_1\rightarrow X_2$
is a non-zero map. For any $t\in [0,1]$,
define maps $q_t:X_1\rightarrow\mathbb R^n$ and
$h_t:X_1\rightarrow X_2$ by
$q_t(x)=tq(x)$, $h_t(x)=\pi\circ\sigma(q_t)$.
For any $t\in[0,1]$, $h_t$ is a continuous
equivariant map from $X_1$ to $X_2$, and $h_t(e)=e$.
We claim that $h_t$ is non-affine for some $t\in(0,1]$.

Suppose this is not the case. Then for each $t\in[0,1]$,
$h_t$ is a continuous homomorphism from $X_1$ to $X_2$.
Let $Y$ denote the set of all continuous
maps from $X_1$ to $X_2$.
Choose any metric $\rho$ on $X_2$ that gives the
topology,
and define a metric $\rho_0$ on
$Y$ by
$$
\rho_0(h_1,h_2)=\sup\{\rho(h_1(x),h_2(x))\mid x\in X_1\}.
$$
The map $t\mapsto h_t$ is continuous with respect to
$\rho_{0}$ and the set of all continuous homomorphisms
from $X_1$ to $X_2$ forms a
discrete subset of $Y$. Hence
$t\mapsto h_t$ is constant, which contradicts the fact that
$h_0=0$ and $h_1\ne0$. This proves that some $h_t$ is
not an affine map.
Since $h_t$ is a
continuous equivariant map
from $\mathsf X_1$ to $\mathsf X_2$
for any $t\in[0,1]$, Theorem~\ref{rigidity}
follows.
\end{proof}

%% \bibliography{refs} %% integrated into single file 2 Feb 2003

\providecommand{\bysame}{\leavevmode\hbox to3em{\hrulefill}\thinspace}

\end{document}